\documentclass[12pt]{amsart}    
\usepackage{latexsym}
\usepackage{amssymb}

\numberwithin{equation}{section} 

\font\tengothic=eufm10 scaled\magstep 1
\font\sevengothic=eufm7 scaled\magstep 1
\newfam\gothicfam
          \textfont\gothicfam=\tengothic
          \scriptfont\gothicfam=\sevengothic

\newcommand{\s}{\; | \;}

\newcommand{\Z}{\mathbb{Z}}

\newcommand {\PP}{\mathbb{P}}
\DeclareMathOperator{\chara}{char}

\DeclareMathOperator{\codim}{codim}
\DeclareMathOperator{\depth}{depth}

\DeclareMathOperator{\soc}{Soc}

\DeclareMathOperator{\Proj}{Proj}

\DeclareMathOperator{\pnt}{\raise 0.5mm \hbox{\large\bf.}}

\newtheorem{theorem}{Theorem}[section]
\newtheorem{lemma}[theorem]{Lemma}
\newtheorem{proposition}[theorem]{Proposition}

\newtheorem{conjecture}[theorem]{Conjecture}

\newtheorem{prop-def}[theorem]{Proposition and Definition}

\theoremstyle{definition}
\newtheorem{remark}[theorem]{Remark}

\newtheorem{rem-def}[theorem]{Definition and Remark}

\title[Gorenstein Hilbert functions in codimension four]
{A characterization of
Gorenstein Hilbert functions in codimension four with small initial degree}

\author[Juan Migliore]{Juan Migliore${}^*$}
\address{
Department of Mathematics, University of Notre Dame, Notre Dame, IN
46556, USA}
\email{Juan.C.Migliore.1@nd.edu}
\author[Uwe Nagel]{Uwe Nagel${}^{+}$}
\address{Department of Mathematics,
University of Kentucky, 715 Patterson Office Tower,
Lexington, KY 40506-0027, USA}
\email{uwenagel@ms.uky.edu}

\author[Fabrizio Zanello]{Fabrizio Zanello${}^\dagger$}
\address{
Department of Mathematical Sciences,
Michigan Technological University,
Houghton, MI 49931-1295, USA}
\email{zanello@math.kth.se}

\thanks{${}^*$ Part of the work for this paper was done while the first
author was sponsored by the National Security Agency under Grant
Number H98230-07-1-0036. \\
${}^+$ Part of the work for this paper was done while the second
author was sponsored by the National Security Agency under Grant
Number H98230-07-1-0065. He also acknowledges partial support from, and the hospitality of, the Institute for Mathematics  \& its Applications at the University of Minnesota. \\
${}^\dagger$ Part of the work  was done when the third author was funded by the G\"oran Gustafsson
Foundation.}

\begin{document}

\begin{abstract}
The main goal of this paper is to characterize the  Hilbert functions
of all (artinian) codimension 4 Gorenstein algebras that have at least two independent relations of degree four. This includes all codimension 4 Gorenstein algebras whose initial relation is of degree at most 3. Our result shows that those Hilbert functions are exactly the so-called
{\em SI-sequences} starting with $(1,4,h_2,h_3,...)$, where $h_4 \leq 33$. In
particular, these Hilbert functions are all unimodal.

We also establish a more general unimodality result, which relies on the values of the Hilbert function not being too big, but is independent of the initial degree.
\end{abstract}

\maketitle

\section{Introduction}

Since the appearance of Bass's paper \cite{bass} ``On the ubiquity
of Gorenstein rings,'' its title has been amply justified. Indeed,
Gorenstein structures are often involved in duality statements and
they have found numerous applications in many areas of mathematics -
including algebraic geometry, combinatorics, and complexity theory
(see, e.g., \cite{Hu}, \cite{P}, \cite{St0}, \cite{A}, \cite{St-80}, \cite{KS}).

Of particular importance is the case when the Gorenstein ring is a
standard graded algebra $A = k[A_1]$ over a field $k$. In this case,
the dimensions of its graded components are basic invariants and of
fundamental interest. They are comprised in the Hilbert function. It
can be described by finitely many pieces of data using the fact that its
generating function, the Hilbert series, is rational and can be
uniquely  written as:
\[
\sum_{j \geq 0} \dim_k A_j \cdot z^j =: \frac{h_0 + h_1 z + \cdots + h_e z^e}{(1-z)^{\dim A}}.
\]
The vector $(h_0,\ldots,h_e)$  consists of  positive integers, and is
called the {\em $h$-vector} of $A$. By duality it is symmetric about
the middle, and in particular $h_e = h_0 = 1$. We call $h_1$ the {\em
codimension} of $A$. A central problem is to put restrictions on, or
even to characterize, the possible Gorenstein $h$-vectors.  This is easy
if $h_1 \leq 2$ because then $A$ is a complete intersection. In
codimension 3 the classification is due to Stanley who showed in
\cite{St1} (see also \cite{Za}) that in this case the $h$-vectors are exactly the
so-called {\em SI-sequences}  with $h_1 = 3$. SI-sequences exist in every
codimension and are defined by a simple numerical condition (see
Section \ref{sec-prel}). It is known (see \cite{Ha}, \cite{MN}, or
\cite{CI}) that every SI-sequence is a Gorenstein $h$-vector.
However, Stanley's Example 4.3 in \cite{St1} shows that the
converse is not true. In fact, the $h$-vector in this example has
codimension 13 and is not even unimodal. Recall that an $h$-vector
is said to be {\em unimodal} if it is never strictly increasing
after a strict decrease. Later on, in each codimension $\geq 5$, examples of
$h$-vectors  were found that are not unimodal (see \cite{BI}, \cite{BL}, \cite{Bo}). Since then the problem has been open whether
non-unimodal Gorenstein $h$-vectors of codimension 4 exist (see,
e.g., \cite{BI}, \cite{IK}, \cite{MV}, \cite[Problem 2.19]{V}, or \cite[p.\ 66]{St2}).

Strengthening the assumption on the ring, the analogous question
becomes meaningful in every codimension. Indeed, Stanley conjectured:

\begin{conjecture}[Stanley \cite{St}] If $A$ is a  standard graded Gorenstein {\rm domain}, then its
$h$-vector is unimodal.
\end{conjecture}

\noindent This conjecture is wide open if the codimension is at least four.

The condition of being an SI-sequence is stronger than that of being unimodal.  Thus an even more interesting question (especially if the answer is affirmative) is whether every codimension 4 Gorenstein $h$-vector is an SI-sequence.  Recently, a partial result was obtained  by  Iarrobino and Srinivasan \cite{IS}. They showed that all Gorenstein  $h$-vectors $(1, 4, h_2,\ldots)$ with $h_2 \leq 7$ are SI-sequences, so in particular unimodal.  Thus their result covers Gorenstein algebras of codimension 4 with initial degree 2 and at least a 3-dimensional family of quadrics in the ideal.

The main result of this paper is a characterization of the Gorenstein codimension 4 $h$-vectors if the ideal  contains enough forms of degree four:

\begin{theorem} \label{thm-intro}
Let $I \subset k[x_1,\ldots,x_4] =:R$ be an artinian Gorenstein ideal
where $k$ is a field of characteristic zero. If $I$ contains at least
two independent quartics, then the $h$-vector of $R/I$ is an
SI-sequence.

Therefore, the $h$-vectors of Gorenstein algebras with $h_1 \leq 4$ and $h_4 \leq
33$ are exactly the corresponding SI-sequences.
\end{theorem}

The last part of the theorem follows because each SI-sequence arises as a Gorenstein $h$-vector.

Notice that Theorem \ref{thm-intro} proves Stanley's conjecture for the case of initial degree 3 or less (and even slightly more), in codimension 4. In fact, we also prove a more general unimodality result, thus adding further evidence to this conjecture.  It roughly says that if there is any value of $s$ for which the $h$-vector in degree $s$ is not too large compared to a specific function of $s$, then the  $h$-vector is unimodal
(see Theorem \ref{more general}).

In order to establish our results we introduce new techniques and also
use the classical method of restricting to a  generic hyperplane, though with
a new twist. In fact, Proposition \ref{prop-key} is a cornerstone of
our approach. Though its statement is algebraic, the proof is
geometric in essence because we translate the statement into a question on the number of conditions imposed by a certain zero-dimensional subscheme on a particular linear system. The proof also uses Bertini's theorem and is the only place where the characteristic zero assumption is used explicitly. After these preparations,  developed in
Section~\ref{sec-prel}, we prove our  unimodality results in Section
3. Using those, we finally establish the  stronger result about
SI-sequences in Section 4.


\section{Preliminary results} \label{sec-prel}

In this section we establish some technical results that we need
later on.

Throughout this note, we will consider standard graded $k$-algebras. Such an
algebra is of the form  $A=R/I$ where $R=k[x_1,\dots ,x_r]$ is the
polynomial ring over the field $k$, $\deg x_i = 1$, and $I$ is a
homogeneous ideal of $R$. Its Hilbert function is $h_A : \Z \to
\Z,\; h_A (j) = \dim_k [A]_j$. To study it, we may and will assume
that the field $k$ is infinite.  In fact, later we will assume that
$k$ has characteristic zero, but this is used explicitly only in the proof of
Proposition \ref{prop-key}. However, we use the conclusion of
Proposition \ref{prop-key} in many of our arguments. It would be
interesting to find a characteristic-free proof of this result.

If $A$ is an artinian algebra, i.e.\ its Krull dimension is zero,
then its Hilbert function has finite
support and is captured in its {\it $h$-vector} $h = h(A)
=(h_0,h_1,\dots ,h_e)$ where $h_i= h_A (i) > 0$ and $e$ is the last
index with this property. The integer $e$ is called the {\it socle
degree} of $A$ (or of the $h$-vector $h(A)$). Moreover, abusing
notation slightly, we call $h_1$ the {\it codimension} of $A$. The $h$-vector $(h_0,h_1,\dots ,h_e)$ is said to be a {\em differentiable O-sequence up to degree $j$} if its first difference through degree $j$, $(h_0, h_1 - h_0, h_2 - h_1,\ldots, h_j - h_{j-1})$, is an O-sequence in the sense of \cite{Ma} or \cite{St1}.

Recall that the artinian $k$-algebra $A$ is Gorenstein if its socle
$\soc (A) =\lbrace a\in A {\ } \mid {\ } a \cdot (x_1,\ldots,x_r) =0\rbrace $ is 1-dimensional,
i.e.\ $\soc (A) \cong k(-e)$. Its $h$-vector is symmetric, i.e.\ $h_i = h_{e-i}$ for all $i$.
It is called an {\em SI-sequence} if in addition
$(h_0, h_1 - h_0, h_2 - h_1,\ldots,
h_{\left \lfloor \frac{e}{2} \right \rfloor}  - h_{\left \lfloor \frac{e}{2} \right \rfloor - 1})$
is an O-sequence.  In other words, the $h$-vector is a  differentiable
O-sequence up to  degree $\lfloor \frac{e}{2} \rfloor$.  We sometimes say that a Gorenstein algebra $A$ has the {\em Stanley-Iarrobino Property} if its Hilbert function is an SI-sequence.

Let $n$ and $i$ be positive integers.
The {\it i-binomial expansion of n} is $$n_{(i)}=\binom{n_i}{i}+{n_{i-1}}{ i-1}+...+\binom{n_j}{ j},$$ where
$n_i>n_{i-1}>...>n_j\geq j\geq 1$. The
$i$-binomial expansion of $n$ is unique (see, e.g., \cite{BH}, Lemma
4.2.6). Hence, we may define
$$
\begin{array}{rcl}
\displaystyle n_{<i>} & = & \displaystyle \binom{n_i-1}{ i}+\binom{n_{i-1}-1}{ i-1}+...+\binom{n_j-1}{
j}, \\ \\
n^{<i>} & = & \displaystyle \binom{n_i +1}{i+1} + \binom{n_{i-1} +1}{i} + \dots + \binom{n_j +1}{j+1}
\end{array}
$$
where we set $\binom{m}{ q}=0$ whenever $m<q$. Note that, for fixed $i$, $n_{<i>}$ and $n^{<i>}$ are both increasing functions in $n$.  We first recall Macaulay's theorem (cf.\ \cite{BH} Theorem 4.2.10):

\begin{theorem}[Macaulay]
The Hilbert function of $A$ satisfies:
\[
h_A(i+1) \leq (h_A(i))^{<i>}.
\]

\end{theorem}

We are ready to state Green's Restriction Theorem.

\begin{theorem}[Green \cite{green}] \label{green thm}
If $L \in R$ is a generic linear form, then the Hilbert function of $A/L A$ satisfies:
$$
h_{A/L A} (i)\leq (h_A (i))_{<i>}.
$$
\end{theorem}

We now present further preparatory results.

\begin{lemma} \label{lem-again-gor}
If $F \in R_d$ is a homogeneous form of degree $d$ that is not in $I$,
then $R/(I : F)$ is Gorenstein with socle degree $e - d$.
\end{lemma}

\begin{proof}
This is well-known and follows easily by computing
 $\soc( R/(I:F)) \cong k(-e+d)$.
\end{proof}

If $\dim_k A_1 \leq 1$, then $A \cong k[x]/(x^e)$. Thus it is harmless to assume $\dim_k A_1 \geq 2$.

\begin{lemma} \label{lem-gcd}
Let $A = R/I$ be a graded artinian algebra with $\dim_k A_1 \geq 2$. Assume that
$I_t$ has a greatest common divisor (GCD), $F$, of degree $d > 0$. Set $B := R/(I : F)$.
Then:
\begin{itemize}
\item[(a)] If $i \leq t$, then $h_B (i-d) = h_A (i) - \left [ \binom{i+r-1}{r-1} - \binom{i-d+r-1}{r-1} \right ]$.

\item[(b)] If the $h$-vector of $B$ is non-decreasing (resp.\ increasing) in degrees $\leq t-d$ then the $h$-vector of $A$ is non-decreasing (resp.\ increasing)  in degrees $\leq t$.

\item[(c)] If the $h$-vector of $B$ is a differentiable O-sequence up to degree $t-d$ then the $h$-vector of $A$ is a differentiable O-sequence up to degree $t$.

\item[(d)] If $A$ is Gorenstein and if the $h$-vector of $B$ is an SI-sequence, then the $h$-vector of $A$
satisfies the SI condition in all degrees $i \leq t$.

\end{itemize}
\end{lemma}

\begin{proof}
Multiplication by $F$ on $A$ induces the exact sequence:
\[
0 \to (R/(I : F)) (-d) \to A \to R/(I, F) \to 0.
\]
Since by assumption $(I, F)_i = (F)_i$ if $i \leq t$, we get for the
Hilbert functions if $i \leq t$:
\begin{equation} \label{eq-SI}
h_A (i) = h_B (i-d) + h_{R/(F)} (i),
\end{equation}
which implies our first and second assertions.

For (c), the assumption on $B$ provides the existence of an
ideal $J \subset S := k[x_1,\ldots,x_{r-1}]$ such that
\[
h_{S/J} (i) =
\left \{
\begin{array}{ll}
\Delta h_B (i) & \hbox{for } i \leq t-d \\
0 & \hbox{for } i > t-d.
\end{array}
\right.
\]
Let $G \in S$ be a form of degree $d$ and consider the exact sequence
\[
0 \to R/J R (-d) \to R/ J G R \to R/G R \to 0.
\]
Its non-trivial modules have positive depth. Thus we conclude for all $i$ that
\begin{eqnarray*}
h_{S/J G} (i) = \Delta h_{R/ J G R} (i) & = & \Delta h_{R/G R} (i)  +  \Delta h_{R/ JR} (i-d)\\
& = &
\Delta h_{R/(F)} (i) +
\left \{
\begin{array}{ll}
\Delta h_B (i-d) & \hbox{for } i-d  \leq t-d; \\
0 & \hbox{for } i-d  > t-d.
\end{array}
\right.
\end{eqnarray*}
Comparing the above calculation with Equation (\ref{eq-SI})  we conclude that $\Delta h_A$ is an $O$-sequence in degrees $i \leq t$, as desired.

For (d), note that if $t > \lfloor \frac{e}{2} \rfloor$, then by symmetry it is enough to show that $h(A)$ is a differentiable O-sequence up to degree $\lfloor \frac{e}{2} \rfloor$. Hence, it suffices to prove that $h(A)$ is a differentiable O-sequence up to degree $j = \min\{t, \lfloor \frac{e}{2} \rfloor \}$. But $B = R/(I:F)$ is again Gorenstein, by Lemma \ref{lem-again-gor}, and it has socle degree $e-d$.  Therefore its $h$-vector $h(B)$ is an O-sequence up to degree $\lfloor \frac{e-d}{2} \rfloor \geq j-d$. Hence we conclude by using part (c).
\end{proof}

If $V$ is a vector space of forms (either a component of $R$ or a
component of a
graded ideal, depending on the context), we say that a {\em generic} element of
$V$ has the property $P$ if there is a Zariski-open subset of $V$ (or of the
projectivization of $V$) with this property.  We say that a {\em generic choice
of $r$ elements of $V$} has the property $P$ if an open subset of the corresponding
product of $r$ copies of $V$ has this property.

The following is a technical, but crucial, result that will be used in
the next section. Its purpose is to study the possible behavior of a vector space of forms when we reduce modulo two generic linear forms to obtain a vector space of forms in a 2-dimensional polynomial ring.

\begin{proposition} \label{prop-key}
Assume that the field $k$ has characteristic zero.
Let $R = k[x_1,x_2,x_3,x_4]$ and let $J = (F,G_1,G_2) \subset R$ be a homogeneous ideal with three minimal generators, where $\deg F = a \geq 2$ and $\deg G_1 = \deg G_2 = b \geq  a$.  Let $L_1, L_2 \in R$ be generic linear forms.  Then $\dim [R/(J,L_1,L_2)]_b = a-1$ if and only if $F,G_1,G_2$ have a GCD of degree $a-1$.  Otherwise $\dim [R/(J,L_1,L_2)]_b = a-2$.
\end{proposition}

\begin{proof}
Let $\bar R = R/(L_1,L_2)$.  For any homogeneous polynomial $A \in R$, let $\bar A$ be the restriction to $\bar R$, and similarly for an ideal $J$.  So the desired dimensions can be expressed as $\dim [R/(J,L_1,L_2)]_b = \dim [\bar R/ \bar J]_b$.  We will use both formulations.

Observe that $\dim [R/(F,L_1,L_2)]_t = a$ for all $t \geq a-1$.  As long as $G_1 \notin (F)$, it is not hard to see that $\dim [R/(F,G_1,L_1,L_2)]_b = a-1$.  So this lemma gives a criterion for $G_2$ to impose an additional condition modulo $L_1,L_2$, in other words a criterion for the condition that $G_2 \notin (F,G_1,L_1,L_2)$ for generic $L_1,L_2$.  In any case it is clear that the indicated dimension is either $a-1$ or $a-2$.

First suppose that $F,G_1,G_2$ have a non-trivial GCD, $D$.  Let $F' = F/D$, $G_i' = G_i /D$, $i = 1,2$, so $F',G_1', G_2'$ have no GCD of positive degree.  By taking general linear combinations, we may assume that no two of these polynomials have a non-trivial GCD.  Note that $\deg G_1' = \deg G_2' = b-d$, and that $\deg F' = 1$ if and only if $d := \deg D = a-1$.

Now assume further that $d = a-1$.  We first consider $\dim [R/(F',G_1',L_1,L_2)]_{t}$.  Since $\deg F' = 1$, $F',G$ form a regular sequence, and $L_1,L_2$ being  generic, we obtain that this dimension is 1 up to degree $t = b-d-1= b-a$, and 0 afterwards.  Hence we have $G_2' \in (F',G_1',L_1,L_2)$.    It follows that $G_2 \in (F,G_1,L_1,L_2)$, so $\dim [R/(J,L_1,L_2)]_b = a-1$ as desired.

Next assume that $F,G_1,G_2$ have a GCD, $D$, of degree $\leq a-2$.  Consider the exact sequence
\[
0 \rightarrow \bar R/(\bar J : \bar D) (-d) \rightarrow \bar R / \bar J \rightarrow \bar R / (\bar J + (\bar D)) \rightarrow 0
\]
and note that this last  module is equal to $\bar R / (\bar D)$.  Note also that $\bar J : \bar D = (\bar F', \bar G_1', \bar G_2')$.  We obtain
\[
a-2 \leq \dim [\bar R / \bar J]_b = \dim [\bar R / (\bar J : \bar D)]_{b-d} + \dim [\bar R/( \bar D)]_b = \dim [\bar R / (\bar F', \bar G_1', \bar G_2')]_{b-d} + d.
\]
Since $\deg \bar F' = a-d$, if we prove that $\dim [\bar R / (\bar F', \bar G_1', \bar G_2')]_{b-d} = a-d-2$, we will have completed the case where $F,G_1,G_2$ have a GCD of degree $\leq a-2$.

So it remains to treat the case where  $F,G_1,G_2$ have no common non-trivial GCD, and  in fact pairwise they have no GCD.
\smallskip

\noindent
{\bf Case 1:}  Assume that $F,G_1,G_2$ form a regular sequence.  Then $(F,G_1,G_2,L_1) = (J,L_1)$ is a general artinian reduction of $J = (F,G_1,G_2)$ and we have an exact sequence
\[
[R/(J,L_1)]_{b-1} \stackrel{\times L_2}{\longrightarrow} [R/(J,L_1)]_b \rightarrow [R/(J,L_1,L_2)]_b \rightarrow 0
\]
By the Weak Lefschetz Property (WLP) for complete intersections of codimension three (cf.\ \cite{HMNW} -- here we need characteristic zero), the map $(\times L_2)$ has maximal rank.  Since up to degree $b-1$ the ideal $J$ has at most one minimal generator, namely $F$, we compute
\[
\begin{array}{rcl}
\dim [R/(J,L_1)]_{b-1} & = & \binom{b-1+2}{2} - \binom{b-1-a+2}{2} \\ \\
\dim [R/(J,L_1)]_b & = &  \binom{b+2}{2} - \binom{b-a+2}{2} -2
\end{array}
\]
and hence
\[
\dim [R/(J,L_1,L_2)]_b = (b+1) - (b-a+1) -2 = a-2
\]
as desired.
\smallskip

\noindent
{\bf Case 2:} Finally assume that the ideal $(F,G_1,G_2)$ has
codimension two but we still assume that pairwise the generators are
regular sequences.  If $J$ is the saturated ideal of an arithmetically
Cohen-Macaulay curve in $\PP^3$ then $L_1,L_2$ is a regular sequence
on $R/J$ and the result follows easily.  So  we may assume that
$\depth R/J = 0$ or $1$.  In any case, though, $J \supset (F,G_1)$,
which defines a complete intersection curve of degree $ab$.  Since
$G_2 \notin (F,G_1)$, the curve defined by $J$ must have degree
strictly less than $ab$.  (This comes by  comparing the
  top-dimensional part of $(F,G_1,G_2)$ with the unmixed ideal
  $(F,G_1)$.)

We want to show that $G_2 \notin (F,G_1,L_1,L_2)$.  Now, $(F,G_1,L_1)/(L_1)$
is the saturated ideal of a complete intersection zero-dimensional
scheme in $\mathbb P^2 = \Proj\ (R/(L_1))$ of degree $ab$.  Furthermore, we claim that
$G_2 \notin (F,G_1,L_1)$.  Indeed, since $L_1$ is generic, the ideal
$(F,G_1,G_2,L_1)$ in large degree agrees with the saturated ideal of a
zero-dimensional scheme of degree $< ab$.


For simplicity, then, we rephrase the problem as follows.  Let
$(F',G_1') \subset k[x,y,z]$ be a complete intersection ideal with
$\deg F' = a, \deg G_1' = b$, and let $L$ be a generic linear form.
Let $G_2'$ be a homogeneous polynomial in $k[x,y,z]_b$ which is not in
$(F',G_1')$.  Then we claim that $G_2' \notin (F',G_1',L)$.   By
semicontinuity, it suffices to produce one $L$ for which this is the
case.

We begin with the  ideal $(F',G_1', G_2')$ and look for $L$ shortly.
This ideal defines a set of $< ab$ points (but it is not necessarily
saturated).  Replacing $G_1'$ and $G_2'$ by sufficiently general
elements of this ideal in degree $b$ (but by abuse of notation still
calling them $G_1'$ and $G_2'$), we see that the pencil of curves in
$\mathbb P^2$ spanned by $G_1'$ and $G_2'$ has base locus consisting
of finitely many points (in fact $b^2$, counted with multiplicities).
To show that $G_2' \notin (F',G_1',L)$ (where $L$ is still to be
introduced) will be equivalent to showing that in $k[x, y, z]$ there
is no relation
\begin{equation}\label{rel}
G_2' = AF' + \lambda G_1' + BL,
\end{equation}
where $\deg A = b-a$, $\deg B = b-1$ and $\lambda$ is a scalar.  We
will suppose otherwise and obtain a contradiction, i.e. we will
suppose that no matter what $L$ we use, there is a relation of this
form.

We may  assume that the field $k$ is algebraically closed.
Now we consider the pencil of curves in $\mathbb P^2$ spanned by
$G_1', G_2'$.  We will abuse notation and refer to curves and defining
polynomials by the same name.  Let $P$ be a point on $F'$ away from
the (finite) base locus of the pencil.  There is a unique element,
$H$, of the pencil vanishing at $P$.  We will show that $ {F'}$ must
divide $H$, thus obtaining a contradiction to the assumption that
$G_2' \notin (F',G_1')$. To this end we distinguish several cases.
\smallskip

{\em Case 2.1}.  Suppose that $F'$ is reduced. We may assume that the
point $P$ has coordinates $(0 : 0 : 1)$. Let $Q = (\alpha : \beta :
\gamma ) \neq P$ be any other point lying on $F'$. Then the linear
form $L = \beta x -  \alpha y$ is not zero, and $P$ as well as $Q$ lie
on $L$.
If, for this $L$, there is no Relation (\ref{rel}), then we are
done. Otherwise, we evaluate this relation at $P$ and conclude that
$H$ is in the saturated ideal $(F', L)$. (Note that this is also true
  if $P$ and $Q$ lie on the same linear component of $F'$ because then
  $(F', L) = (L)$).  It follows that $Q$ lies on $H$. Since every
point $Q$ of $F'$ lies on $H$ and $F'$ is reduced, this implies that
$F'$ divides $H$, which is the desired contradiction.
\smallskip

{\em Case 2.2}.  Now suppose that $F'$ is not reduced, but still has
at least one non-linear component. Let $P$ again be a point on $F'$,
but not lying on any linear component.  Consider an irreducible
factor, $g$, of $F'$ and denote by $e_0 \geq 1$ its multiplicity. Let
$L$ be the line through $P$ and any point $Q \neq P$ lying on
$g$. Then the intersection
multiplicity of  $L$ with $F'$ at $Q$ is at
least $e_0$. Remembering that $H$ is determined only by $P$, the
argument in Case 2.1 shows that $H$ vanishes with multiplicity at
least $e_0$ at each such point $Q$. It follows that $g^{e_0}$ divides $H$,
thus so does $F'$.
\smallskip

{\em Case 2.3}.  If the support of $F'$ is a union of (at least two)
lines,  we choose $P$ such that it lies on exactly one of the
components. Then
the above argument applies to all components of $F'$ but the line
through $P$. Let $P'$ be a point on $F'$ that is not on this line. Then
$P'$ also lies on the curve $H$ determined by $P$. Hence the pencil of
lines through $P'$ can be used to fill in the missing component, and
we again get that $F'$ is contained in $H$.\smallskip

{\em Case 2.4}.  Finally, suppose that $F'$ is supported on a line.
Let $P$ be any point on this line away from the base locus of the
pencil (which, by our assumption, is a finite set of points, so the
supporting line is not part of the base locus).  There is a unique
element of the pencil, $H_P$, containing $P$.  By assumption, if $L$
is a generic line through $P$ then $H_P$ contains the degree $a$
intersection of $L$ with $F'$ (which is supported at $P$).  But this
is true for infinitely many such $L$ containing $P$, and $H_P$ remains
unchanged, which means that $H_P$ has a singular point of multiplicity
at least $a$ at $P$.  There are infinitely many choices for $P$ on
$L$, and the same conclusion holds for each.  This can mean one of two
things.  If $H_P$ remains constant as $P$ varies on $L$ then $F'$
divides $H_P$ and we have a contradiction as above.  If $H_P$ varies,
then the general element of the pencil is singular away from the base
locus, contradicting Bertini's theorem (see, e.g., \cite{kleiman}).

This completes the proof.
\end{proof}

\begin{remark}
  \label{rem-ext-2.5}
(i) The characteristic zero assumption on the field is used in two places: to apply Bertini's theorem and to use the result in \cite{HMNW} that each complete intersection of codimension three has the Weak Lefschetz Property. In fact, the latter statement is false in every positive characteristic $p$. Indeed, let $A = k[x_1, x_2, x_3]/(x_1^p, x_2^p, x_3^p)$ with $\chara k = p > 0$ and consider the multiplication by a generic linear form $L$ on $A$
\[
[A]_{i-1} \stackrel{\times L}{\longrightarrow} [A]_i.
\]
It does not have maximal rank if $i = p$ because then the residue class of $L^{p-1}$ is a non-trivial element in the kernel. If $A$ were to have the Weak Lefschetz Property, then this map would be injective if $j \leq \frac{3p-3}{2}$ and surjective otherwise. However, one can even show that the map does not have maximal rank for all $i = p, p+1, \ldots, 2p-2$.

(ii)
It would be interesting to extend Proposition \ref{prop-key} to ideals $J$ with more than 3 generators. More precisely, assume $J = (F, G_1,\ldots,G_m)$ has $m+1$ minimal generators and $\deg G_i = b \geq a = \deg F$. If the forms $G_1,\ldots,G_m$ are generic, then we have
\[
\dim [R/(J, L_1, L_2)]_b = \max \{0, a - m\},
\]
where $L_1, L_2$ are again generic linear forms. The problem is to find mild conditions on the forms $F, G_1,\ldots,G_m$ that imply the same conclusion. This is easy if $m=1$ and Proposition \ref{prop-key} gives the complete solution if $m=2$.
\end{remark}


\section{Unimodality}

We begin this section with a result that provides unimodality for artinian Gorenstein algebras with ``small'' initial degree. In order to use Proposition \ref{prop-key}, we make throughout the remainder of this note the assumption that the base field $k$ has characteristic zero.

\begin{theorem} \label{unimodal thm}
Let $R/I$ be an artinian Gorenstein algebra with $h$-vector $h = (1,h_1,\dots,h_e)$ and $h_1 \leq 4$.  If  $h_4 \leq 33$ then $h$ is unimodal.
\end{theorem}

\begin{proof}
 We will use induction on the socle degree, $e$, of $R/I$.  Certainly the claim is true when $e=1$.  Now assume to the contrary that $R/I$ has the  least socle degree $e \geq 2$ among  Gorenstein algebras whose  $h$-vectors satisfy our  assumption and are not unimodal.
Let $a := \min \{j \in \Z  \s I_j \neq 0\}$ be the initial degree of $I$.  Our assumption says, in particular, that $a \leq 4$.  By Stanley's Theorem 4.2 in \cite{St1}, the claim is true  if $h_1 \leq 3$. Thus, we may
 assume  $h_1 = 4$ and $2 \leq a \leq 4$. Let $F \in I$ be a minimal generator of degree
 $a$.  Note that if $h_2 = 9$ or $h_3 = 19$  then $F$ is unique up to scalar multiple.

Let $L_1$ be a generic linear form and let $b = (1,b_2,b_3,\dots)$ be the  $h$-vector
 of the Gorenstein algebra $R/(I:L_1)$ (note that for the sake of the notation below, for
 each $j$ we write $b_j$ for the Hilbert function in degree $j-1$).
 Let $c = (1,3,c_2,\dots)$ be the $h$-vector of $R/(I,L_1)$.    We have the computation:

\begin{center}

\begin{equation} \label{hbc diag}
 \begin{tabular}{c|ccccccccccccccccccccccccccccccccccccccccc}
 {deg} & 0 & 1 & 2 & 3 & 4  & \dots &  $\left \lfloor \frac{e}{2} \right \rfloor -1$ & $\left \lfloor \frac{e}{2} \right \rfloor$ & $\left \lfloor \frac{e}{2} \right \rfloor +1$ &  $\left \lfloor \frac{e}{2} \right \rfloor +2$ & \dots & $e-2$ & $e-1$ & $e$ \\[2pt] \hline
 & 1 & 4 & $h_2$ & $h_3$ & $h_4$ & \dots & $h_{\left \lfloor \frac{e}{2} \right \rfloor -1}$ & $h_{\left \lfloor \frac{e}{2} \right \rfloor}$ & $h_{\left \lfloor \frac{e}{2} \right \rfloor +1}$ & $h_{\left \lfloor \frac{e}{2} \right \rfloor +2}$ & \dots & $h_{e-2}$ & 4 & 1  \\
  & $-$ & 1 & $b_2$ & $b_3$ & $b_4$  & \dots & $b_{\left \lfloor \frac{e}{2} \right \rfloor-1}$ & $b_{\left \lfloor \frac{e}{2} \right \rfloor}$ & $b_{\left \lfloor \frac{e}{2} \right \rfloor+1}$ & $b_{\left \lfloor \frac{e}{2} \right \rfloor +2}$ & \dots & $b_{e-2}$ & $b_{e-1}$ & 1    \\[2pt] \hline
 & 1 & 3 & $c_2$ & $c_3$ & $c_4$ & \dots & $c_{\left \lfloor \frac{e}{2} \right \rfloor -1}$ & $c_{\left \lfloor \frac{e}{2} \right \rfloor}$ & $c_{\left \lfloor \frac{e}{2} \right \rfloor +1}$ & $c_{\left \lfloor \frac{e}{2} \right \rfloor +2}$ &
 \dots
\end{tabular}
 \end{equation}
 \end{center}

\noindent  We want to stress the following:

\begin{itemize}

\item The first and second rows are symmetric; specifically, we have  $h_i = h_{e-i}$ and  $b_i = b_{e+1-i}$ for all $i$.

\item For all $i$ we have both $b_i \leq h_{i-1}$ and $b_i \leq h_i$.  Indeed, this follows from the exact sequence
\begin{equation*} 
\begin{array}{ccccccccccccccccccccccc}
0 & \rightarrow & [( {I:L_1}) /{I}](-1) & \rightarrow & [R/I](-1) & \stackrel{\times L_1}{\longrightarrow} & R/I & \rightarrow & R/(I,L_1) & \rightarrow & 0 \\
&&&& \hfill \searrow && \nearrow \hfill \\
&&&&& [R/(I:L_1)](-1) \\
&&&& \hfill \nearrow && \searrow \hfill \\
&&&& 0 &&  \hfill 0
\end{array}
\end{equation*}

\item For all $i$, we have $c_i = h_i - b_i$.

\item As a result of the last observation, $(I:L_1)$ has initial degree less than or equal to that of $I$, and has socle degree one less, as noted in Lemma \ref{lem-again-gor}.  Thus, the induction  hypothesis applies to $I : L_1$. Furthermore, if $I : L_1$ has initial degree 1, then, as noted above,  Stanley's theorem \cite{St1} gives that $R/(I: L_1)$ has unimodal $h$-vector.

\end{itemize}

Since  $h$ is not unimodal, there is a
least integer $i \leq \frac{e}{2} - 1$ such that $h_i > h_{i+1}$. Note that we must have $i + 1 \geq a$, thus $e \geq 2a$.  Moreover, by symmetry,
$h_i > h_{i+1}$ is equivalent to $h_{e-1-i} < h_{e-i}$.    Since by induction $(b_1,\ldots, b_e)$
is unimodal, and since $e-i \geq \frac{e}{2} +1$, we have in particular that $b_{e-i-1} \geq b_{e-i}$, so

\begin{equation} \label{eq1}
c_{e-i-1} = h_{e-i-1} - b_{e-i-1} < h_{e-i} - b_{e-i} = c_{e-i}.
\end{equation}


We now consider the algebra $R/(I,L_1)$ and another generic linear form $L_2$.  Even though $R/(I,L_1)$ is
not Gorenstein, we can still make the same computation:

\begin{center}

\begin{equation} \label{cdf diag}
 \begin{tabular}{c|ccccccccccccccccccccccccccccccccccccccccc}
 {deg} & 0 & 1 & 2 & 3 & 4 & 5 & \dots &   $\left \lfloor \frac{e}{2} \right \rfloor -1$ & $\left \lfloor \frac{e}{2} \right \rfloor$ & $\left \lfloor \frac{e}{2} \right \rfloor +1$ &  $\left \lfloor \frac{e}{2} \right \rfloor +2$ & \dots \\[2pt] \hline
& 1 & 3 & $c_2$ & $c_3$ & $c_4$ & $c_5$ & \dots & $c_{\left \lfloor \frac{e}{2} \right \rfloor -1}$ & $c_{\left \lfloor \frac{e}{2} \right \rfloor}$ & $c_{\left \lfloor \frac{e}{2} \right \rfloor +1}$ & $c_{\left \lfloor \frac{e}{2} \right \rfloor +2}$ & \dots  \\
  & $-$ & 1 & $d_2$ & $d_3$ & $d_4$ & $d_5$ & \dots & $d_{\left \lfloor \frac{e}{2} \right \rfloor-1}$ & $d_{\left \lfloor \frac{e}{2} \right \rfloor}$ & $d_{\left \lfloor \frac{e}{2} \right \rfloor+1}$ & $d_{\left \lfloor \frac{e}{2} \right \rfloor +2}$   \\[2pt] \hline
 & 1 & 2 & $f_2$ & $f_3$ & $f_4$ & $f_5$ & \dots & $f_{\left \lfloor \frac{e}{2} \right \rfloor -1}$ & $f_{\left \lfloor \frac{e}{2} \right \rfloor}$ & $f_{\left \lfloor \frac{e}{2} \right \rfloor +1}$ & $f_{\left \lfloor \frac{e}{2} \right \rfloor +2}$ & \dots
 \end{tabular}
 \end{equation}
 \end{center}

\noindent Here the $d$ row corresponds to $R/((I,L_1) : L_2)$ and
the $f$ row to $R/(I,L_1,L_2)$. We  observe that the fact that
 $d_{e-i} \leq c_{e-1-i}$, together with Inequality (\ref{eq1}), implies
\begin{equation} \label{eq-f}
f_{e-i} > 0.
\end{equation}

\noindent We note again that $e-i \geq \frac{e}{2}+1$.


Now we study the bottom line of the above table.  Because it represents the Hilbert function of an artinian quotient of $R/(L_1,L_2) \cong k[x,y]$, we have the following observations.

\begin{enumerate}
\item For all $j$, if $I$ has only one minimal generator, of degree $a$, in degrees $\leq j$ then $f_j= a$.

\item Past degree $a$, $\{ f_j \}$ is non-increasing, and in particular $f_{\left \lfloor \frac{e}{2} \right \rfloor+1} \geq f_{e-i}  > 0$.

\item Thanks to Green's theorem (Theorem \ref{green thm}), if $I$ has at least two minimal generators in degrees $\leq j$ then $f_j \leq  a-1$.  Furthermore,  by Proposition \ref{prop-key}, $I_j$ has a GCD of degree $a-1$ if and only if $f_j = a-1$.

\end{enumerate}

\noindent What, then, may be the sequence $\{ f_j \}$? We distinguish three cases according to the initial degree $a$.

\bigskip

\noindent {\bf Case 1:} $a = 2$.

Certainly, as long as $I$ is principal its $h$-vector is a differentiable O-sequence, hence unimodal.  At some point $j \leq \frac{e}{2}$ there is a second generator, causing $f_j \leq 1$.  Since $f_{\left \lfloor \frac{e}{2} \right \rfloor +1} = 1$,  the sequence $\{ f_j \}$ is of the form $(1,2,2,\dots,2,1,1,\dots,1)$.  Since, by assumption of non-unimodality, the second generator comes in degree $\leq \frac{e}{2}$, we have the values $(f_{\left \lfloor \frac{e}{2} \right \rfloor}, f_{\left \lfloor \frac{e}{2} \right \rfloor+1}) = (1,1)$.  This represents maximal growth of the Hilbert function of the  algebra $R/(I,L_1,L_2)$ from degree $\left \lfloor \frac{e}{2} \right \rfloor$ to degree $\left \lfloor \frac{e}{2} \right \rfloor+1$.  Thanks to Davis's results (see \cite{davis} or \cite{BGM}), it follows that $(I,L_1,L_2)/(L_1, L_2)$ has a GCD of degree 1 in degree $\left \lfloor \frac{e}{2} \right \rfloor+1$.
We claim that by the genericity of $L_1$ and $L_2$, this implies that $I$ also has a GCD of degree 1 in the same
degree.  Indeed, the latter statement means precisely that the ideal,
$J$, generated by all forms of degree $\lfloor \frac{e}{2} \rfloor +1$ in
$I$ has codimension 1 and degree 1.  If the assertion were false and
the codimension of $J$ were at least 2 then, by an elementary version of Bertini's theorem that is valid in arbitrary characteristic, $\codim (J,L_1,L_2)/(L_1, L_2) \geq 2$, a contradiction.  But since $J$ then has codimension 1, the degree is
preserved when cutting with two generic hyperplanes, and the
assertion is proved.  We use this type of argument without further comment in several places later on to lift information from $(I, L_1, L_2)/(L_1, L_2)$ to $I$.

Now we apply Lemma \ref{lem-gcd} (b).  Let $t = \left \lfloor \frac{e}{2} \right \rfloor +1$.  We have that $I_t$ has a GCD, $D$, of degree 1 and the $h$-vector of $B = R/(I:D)$ is non-decreasing in degrees $\leq t-1$ (since it is Gorenstein of socle degree $e-1$, using the induction hypothesis).  Hence the $h$-vector of $R/I$ is unimodal, contradicting our assumption.
\bigskip

\noindent {\bf Case 2:} $a = 3$.

Now the sequence $\{ f_j \}$ begins $(1,2,3,f_3 ,\dots)$ and $f_{e-i} \geq 1$ (recall $e-i \geq \frac{e}{2} +1$).  If we have a GCD in any degree, it has degree $\leq 2$.  We have the following possibilities for the pair $(f_{e-i-1}, f_{e-i})$:

\begin{itemize}
\item $(3,3)$, $(3,2)$ or $(3,1)$. In any of these cases the second minimal generator for $I$ comes in degree $\geq \frac{e}{2}+1$, which is impossible if $h$ is not unimodal.

\item $(2,2)$ or $(1,1)$.  Since $\frac{e}{2} \geq a = 3$, this again represents maximal growth of the Hilbert function of $R/(I,L_1,L_2)$, thus $I$ has a GCD of degree 2 or 1, respectively, in degree $e-i$, and so also in degree $\left \lfloor \frac{e}{2} \right \rfloor+1$.  Hence the same argument as in the case $a=2$ gives unimodality.

\item $(2,1)$.  We consider two subcases:

\begin{enumerate}
\item $e-i > \frac{e}{2} +1$.   In this case the sequence $(f_{e-i-2}, f_{e-i-1}, f_{e-i})$ is either $(3,2,1)$ or $(2,2,1)$.  In the first case, the second minimal generator of $I$ comes in degree $\geq \frac{e}{2} +1$, so unimodality follows immediately.  In the latter case $I$ has a GCD of degree 2 in degree $e-i-1 > \frac{e}{2}$, so the same argument as above applies.

\item $e-i = \frac{e}{2} +1$. Note that
necessarily $e$ is even in this case. If $I$ has only two minimal generators up to degree $\frac{e}{2}$, then clearly $R/I$ is unimodal. Otherwise, we claim that $I_{\frac{e}{2}}$ has a GCD, $D$, of degree 2. In fact, if $f_{\frac{e}{2}-1} = 2$, then this represents maximal growth from $f_{\frac{e}{2}-1}$ to $f_{\frac{e}{2}}$, and the claim follows. The only other possibility is $f_{\frac{e}{2}-1} = 3.$ Then
we apply Proposition \ref{prop-key} to the ideal $J$ that is generated by the quadratic minimal generator of $I$ and two generic forms of degree $\frac{e}{2}$ in $I$. It shows that $J_{\frac{e}{2}}$ has a GCD of degree 2. But by the generic choice of the degree $\frac{e}{2}$ generators of $J$, it follows that in fact $I_{\frac{e}{2}}$ has a GCD of degree 2, as claimed.  Note that by induction, $R/(I:D)$ is unimodal and so in particular has an $h$-vector that is non-decreasing up to degree $\frac{e}{2}-2$ (even $\frac{e}{2}-1$).  Hence Lemma \ref{lem-gcd} gives that  $R/I$ has an $h$-vector that is non-decreasing up to degree $\frac{e}{2}$, and thus is unimodal.
\end{enumerate}
\end{itemize}

\noindent {\bf Case 3:} $a = 4$.

Note that in this case we have $h_3 = 20$ and $h_4 \leq 33$.  Furthermore, note that Green's theorem, applied (twice) to degree 4, shows that $f_4 \leq 3$, so it follows that $f_j \leq 3$ for $j \geq 4$.  Observe also that $f_j = 4$ if and only if $j = 3$.
In order to reduce the problem, we first assume that $e-i > \frac{e}{2} +1$. Since $\frac{e}{2} \geq a = 4$, we then get that $f_{e- i - 2} \leq 3$.
   Thus, the possible values for the tuple $(f_{e-i-2}, f_{e-i-1}, f_{e-i})$ now are

\begin{itemize}
\item $(3,3,3)$, $(3,3,2)$, $(3,3,1)$, $(3,2,2)$, $(3,1,1)$, $(2,2,2)$, $(2,2,1)$, $(2,1,1)$, $(1,1,1)$.  In all these cases, there is a GCD in degree $e- i - 1 > \frac{e}{2}$, so Lemma \ref{lem-gcd} shows that $h_{R/I}$ is unimodal.

\item $(3,2,1)$.  In this case, Proposition \ref{prop-key} shows that either $I_{e-i-1}$ has a GCD of degree 3 or else $I$ has only two minimal generators up to degree $e-i-1 > \frac{e}{2}$.  As we have seen, either of these forces $h_{R/I}$ to be unimodal.
\end{itemize}

We conclude: {\em  If $h_{R/I}$ is non-unimodal with minimal possible socle degree, $e$, among artinian Gorenstein algebras with non-unimodal Hilbert function and ideals with initial degree 4, then $e$ is even and the non-unimodality occurs precisely in the middle: $h_{\frac{e}{2} -1} > h_{\frac{e}{2}} < h_{\frac{e}{2} +1}$. }

Returning to our analysis of the triples, which are now $(f_{\frac{e}{2}-1}, f_{\frac{e}{2}}, f_{\frac{e}{2}+1})$, the same sort of arguments as above eliminates situations where the second entry is 3, or where two entries are equal.  The only cases that we have to deal with are as follows:

\begin{itemize}
\item $(4,2,1)$.  As noted above, the 4 in the first entry implies that $\frac{e}{2} -1 = 3$, i.e. $e = 8$.  Hence the Hilbert function of $R/I$ is $(1,4,10,20,h_4, 20,10,4,1)$.  Non-unimodality forces $h_4 < 20$.  But then Green's theorem implies that $f_4 = f_{\frac{e}{2} -1}  \leq 1$, a contradiction.

\item $(3,2,1)$.  The remainder of the proof will be devoted to eliminating this case.

\end{itemize}

We observe that
\[
\begin{array}{rcl}
1 & = & f_{\frac{e}{2} +1}  \\ [2pt]
& = & c_{\frac{e}{2} +1} - d_{\frac{e}{2} +1} \\ [2pt]
& \geq & c_{\frac{e}{2} +1} - c_{\frac{e}{2}} \\ [2pt]
& = & (h_{\frac{e}{2} +1} - h_{\frac{e}{2}} ) + ( b_{\frac{e}{2}} - b_{\frac{e}{2} +1} ) \\ [2pt]
& \geq & h_{\frac{e}{2} +1} - h_{\frac{e}{2}} \\ [2pt]
& > & 0.
\end{array}
\]
It follows that {\em  if $h_{R/I}$ is non-unimodal with minimal possible socle degree, $e$, among artinian Gorenstein algebras with non-unimodal Hilbert function, then $e$ is even and $h_{\frac{e}{2} -1} = h_{\frac{e}{2} +1} = h_{\frac{e}{2}} +1$ as well as $b_{\frac{e}{2}} = b_{\frac{e}{2} + 1}$ and $d_{\frac{e}{2} + 1} = c_{\frac{e}{2}}$. }

Now, we will combine (\ref{hbc diag}) and (\ref{cdf diag}), and use the above considerations to simplify part of the new diagram.  Let $m = d_{\frac{e}{2}}$, $n = b_{\frac{e}{2}}$ and $h = h_{\frac{e}{2}}$.  Then we have the diagram

\begin{center}

\begin{equation} \label{new diag}
 \begin{tabular}{c|ccccccccccccccccccccccccccccccccccccccccc}
 {deg} &  & \dots &  $\frac{e}{2} -1$ & $\frac{e}{2}$ & $\frac{e}{2} +1$ &  $\frac{e}{2} +2$ & \dots \\[3pt] \hline
 & & \dots & $h+1$ & $h$ & $h+1$ & $h_{\frac{e}{2} +2}$ & \dots  \\ [3pt]
  &  & \dots & $(\leq n)$ & $n$ & $n$ & $(\leq n)$ & \dots   \\[3pt] \hline
 & & \dots & $(\geq m+3)$ & $m+2$ & $m+3$ & $c_{\frac{e}{2} +2}$ &
 \dots \\ [3pt]
 && \dots & & $m$ & $m+2$ & $c_{\frac{e}{2} +2}$ \\ [3pt] \hline
 && \dots & 3 & 2 & 1 & 0
 \end{tabular}
 \end{equation}
 \end{center}

 If $h_4 \leq 30$ then $f_4 \leq 2$ (by Green's theorem).  This means $\frac{e}{2} = 4$ and we have a contradiction as above.  So we have $h_4 = 31$, 32 or 33 and $e \geq 10$.  But in fact, if $e = 10$ then the Hilbert function of $R/I$ is $(1,4,10,20,h_4,h_5,h_6,20,10,4,1)$.  Our assumption that $h_4 = 31$,  32 or 33, plus non-unimodality, means $h_5 \leq 32$.  But applying Green's theorem we get $f_5 \leq 1$, again giving a contradiction.  Hence $e \geq 12$.

 It follows that the $f$-line (the Hilbert function of $R/(I,L_1,L_2)$) is $(1,2,3,4,3,3, \dots$ $\dots , 3,2,1)$, where the last 2 occurs in degree $\frac{e}{2} \geq 6$.  Thus $I$ has a GCD, $A$, of degree 3 in all degrees from 4 to $\frac{e}{2} -1$.  In particular,  we get that whenever $4 \leq j \leq \frac{e}{2} - 1$:
\[
h_{R/(I,A,L_1,L_2)} (j) = h_{R/(A,L_1,L_2)} = 3.
\]
Comparing with the Hilbert function of $R/(I,L_1,L_2)$, we conclude that in these degrees
\[
(I, L_1, L_2)_j = (I, A, L_1, L_2)_j.
\]
Since $A$ has degree 3, we get that, for all $j \geq 4$,
\[
(I, L_1, L_2)_j = (I, A, L_1, L_2)_j.
\]
Now applying Proposition \ref{prop-key} using $A$ and two generic forms in $(I,A)$ of degree $\frac{e}{2}$ and lifting the information from $(I, A, L_1, L_2)/(L_1, L_2)$ to $(I, A)$ as in Case 1 of the proof of Theorem \ref{unimodal thm}, we see that $(I,A)$ either has a GCD of degree 2 in degree $\frac{e}{2}$ or else it has only two minimal generators up to degree $\frac{e}{2}$. In the former case we get in particular that $I$ has a GCD of degree 2 in degree $\frac{e}{2}$, hence the above arguments (using Lemma \ref{lem-gcd}) show that $h_{R/I}$ is unimodal. In the latter case, we obtain that the Hilbert function of $R/(I, A)$ is non-decreasing up to degree $\frac{e}{2}$. By our induction hypothesis, $R/(I : A)$ has a unimodal $h$-vector. Given the relation between the Hilbert functions, this means that the Hilbert function of $R/I$ is non-decreasing up to degree $\frac{e}{2}$, and hence it is unimodal  by symmetry.
\end{proof}

 \begin{remark}

 \begin{enumerate}

 \item Theorem \ref{unimodal thm}  shows that Gorenstein ideals with initial degree $\leq 3$ have  quotient rings with unimodal $h$-vector, and nearly proves the same for initial degree 4; the only missing case is $h_4 = 34$.

\item Even assuming that the initial degree is 2, the unimodality conclusion of Theorem \ref{unimodal thm} is false in every codimension five or more.  Indeed, Example 2 of \cite{BI} has codimension 5 and  initial degree 2, and
even the much earlier Example 4.3 of \cite{St1} has codimension 13 and
initial degree 2. Moreover, from the non-unimodal Gorenstein $h$-vector of codimension 5 and initial
degree 2, it is easy to construct (for instance by using inverse systems) non-unimodal
Gorenstein $h$-vectors of any codimension $r\geq 6$ and initial degree 2.

 \end{enumerate}
 \end{remark}

 The method of Theorem \ref{unimodal thm} can be extended to arbitrary initial degree, with some restrictions on the Hilbert function.  We stress that the following theorem does {\em not} immediately include the case $a = 4$ above, because of the restriction $s +1 < \frac{e}{2}$.  We were able to avoid this restriction in Theorem \ref{unimodal thm} by using Proposition \ref{prop-key} {\em twice}, but this relied on the initial degree being 4.  But in any case note that when $s =a = 4$, the condition $h_s \leq 2s^2 +1$ in the following theorem gives exactly $h_4 \leq 33$.

 \begin{theorem}\label{more general}
 Let $R/I$ be an artinian Gorenstein algebra of socle degree $e$ with $h$-vector $h = (1,h_1,\dots,h_e)$.  Assume that there is some positive integer $s$ satisfying $s+1 < \frac{e}{2}$ and $h_s \leq 2s^2 +1$.  Then $R/I$ has a unimodal $h$-vector.
 \end{theorem}

 \begin{proof}
The proof follows the method used in Theorem \ref{unimodal thm},
especially the case $a=4$.  By Theorem \ref{unimodal thm} we may
assume $s \geq 5$.
Observe that
  \[
  2s^2+1 = \binom{s+2}{s} + \binom{s+1}{s-1} + \binom{s}{s-2} +
\binom{s-1}{s-3} -1.
  \]
  By Green's theorem (Theorem \ref{green thm}), the condition that
$h_s \leq 2s^2+1$ says that when we restrict to $R/(L_1,L_2)$ we
obtain $f_s \leq 3$.  In this case we have no information about $h_j$
for $j\leq s-1$.  However, the first part of Theorem \ref{unimodal
thm} made no use of any assumption about the initial degree, so we
still have that if the Hilbert function of $R/I$ fails to be
unimodal in degree $e-i \geq \frac{e}{2}+1$ then
$f_{e-i} > 0$ (cf.\ (\ref{eq-f})).  By the pigeonhole principle, this
forces $1 \leq f_j = f_{j+1} \leq 3$ for some $j$ satisfying $5 \leq
s \leq e-i-3 \leq j \leq e-i-1$. Thus $I_{j+1}$ has a non-trivial
GCD, $D$, of degree $d \in \{1, 2, 3\}$.  We now proceed by induction.
Suppose that $R/I$ fails to be unimodal, and that its socle degree,
$e$, is the smallest among the socle degrees of the non-unimodal
Hilbert functions for which an integer $s$ exists as in
the statement of the theorem.
\medskip

\noindent {\bf Case 1:} Suppose $j+1 \geq \frac{e}{2}$.  Observe that
$R/(I:D)$ is a Gorenstein algebra of socle degree $e-d$, and that
$\dim_k I_s = \dim_k (I:D)_{s-d}$ because $D$ is a common divisor of $I_s$ and $s \geq d$.  Since
$\binom{s+3}{3} - \dim I_s \leq 2s^2+1$, we have
\[
h_{R/(I:D)}(s-d) \leq \binom{s-d+3}{3} - \binom{s+3}{3} + (2s^2+1).
\]
A simple calculation then gives
\[
h_{R/(I:D)}(s-d) \leq
\left \{
\begin{array}{ll}
\frac{3}{2}s(s-1) & \hbox{if $d=1$;} \\ \\
s(s-2) & \hbox{if $d = 2$}; \\ \\
\frac{1}{2}s(s-3) &\hbox{if $d = 3$}.
\end{array}
\right.
\]
Then one can check that $h_{R/(I:D)}(s-d) \leq 2(s-d)^2 +1$ for $s
\geq 5$ (and in fact even $s \geq 3$).  Thus by induction we may
assume that $R/(I:D)$ has a unimodal Hilbert function, and so by Lemma
\ref{lem-gcd} we are done.
\medskip

\noindent {\bf Case 2:} Since $e-i-1 \geq \frac{e}{2}$, the only
remaining possibility is $j+1 = e-i-2 = \frac{e}{2} - 1$ and $(f_{e-i-2}, f_{e-i-1},
f_{e-i}) = (3,2,1)$, where $I_{e-i-2}$ has a GCD of degree 3. This is handled exactly as in the proof
of  Theorem \ref{unimodal thm}.
\end{proof}

\begin{remark}
We hope that our methods will contribute substantially
to the eventual solution of the problem of whether all Gorenstein
algebras of codimension four have unimodal $h$-vectors.  If, as we
hope, it turns out that they do have unimodal $h$-vectors, then an
attempt to prove this based on our approach seems to require a
suitable extension of Proposition \ref{prop-key} (cf.\ Remark \ref{rem-ext-2.5}(ii)).  If, on
the other hand, there is a counterexample, our methods and results
may help guide the way to the right approach to finding it.
\end{remark}


\section{The Stanley-Iarrobino property}

We now show that a similar analysis can be used to show that all codimension 4 Gorenstein $h$-vectors with initial degree $\leq 3$ or $h_4 \leq 33$ are SI-sequences.  Note that we cannot merge this proof with the preceding one because this proof assumes (thanks to Theorem \ref{unimodal thm}) that the sequences are unimodal to begin with. The result is equivalent to Theorem \ref{thm-intro} in the introduction.

\begin{theorem} \label{SI thm}
Let $R/I$ be a Gorenstein algebra  with $h$-vector $h = (1,h_1,\dots,h_e)$ and $h_1 \leq 4$.  If  $h_4 \leq 33$,  then $h$ is an SI-sequence.
\end{theorem}

\begin{proof}
We may assume, by Theorem \ref{unimodal thm}, that $h$ is unimodal.  Furthermore, as before, the case $h_1 \leq 3$ is known, so we may assume that $h_1 = 4$.  Our strategy will be very similar to that of Theorem \ref{unimodal thm}: we restrict (separately) by two generic linear forms $L_1$ and $L_2$, and make computations similar to those above (see (\ref{hbc diag}) and (\ref{cdf diag})).  We show that the failure of $h$ to be an SI-sequence again forces a non-zero entry in the $f$ line, and work with that.

Suppose that $h$ has socle degree $e$. We use the notation of the proof of Theorem \ref{unimodal thm}.
 Assume that $h$  is non-SI in degree $t+2$.  Since $h$ is unimodal, this means that $a + 1 \leq t+2 \leq \frac{e}{2}$ and $h_{i-1} \leq h_{i}$ for $i \leq \frac{e}{2}$, and that the growth from $(h_{t+1} - h_t)$ to $(h_{t+2} - h_{t+1}$) exceeds maximal growth.  In particular, we have
 \begin{equation} \label{h's}
 h_{t+1} - h_t  < h_{t+2} - h_{t+1}.
 \end{equation}

 We make the following general observation:
 \begin{equation} \label{c's}
 \hbox{If $c_i > c_{i-1}$ for any $i$ then $f_i > 0$.}
 \end{equation}

 \noindent Indeed, this follows immediately from the fact that $d_i \leq c_{i-1}$ and $f_i = c_i - d_i$.  We now make some calculations, which we will analyze shortly.

 \medskip

 \noindent {\bf Case 1:}  Assume $c_{e-t} > c_{e-t-1}$.  In this case, by
(\ref{c's}) we have that $f_{e-t} > 0$.  Since $e-t \geq \frac{e}{2} +2 \geq
t+4$, in particular we have $f_{t+4} > 0$. We first claim that $f_{t+2} \leq
3$.   In the proof of Theorem \ref{unimodal thm} we have seen that our
assumption $h_4 \leq 33$ provides $f_4 \leq 3$. Hence $f_{t+2} = 4$ implies
$t+2 = 3$, so $t=1$.  But we have $a+1 \leq t+2$,which in turn gives $a \leq
2$.  However, $f_3 = 4$ implies that $a \geq 4$,  a contradiction. We conclude
that $f_{t+2} \leq 3$, as claimed. Thus we have the following possibilities for
$(f_{t+2}, f_{t+3}, f_{t+4})$:
\begin{eqnarray*}
&& (3,3,3), \  (3,3,2) , \ (3,3,1),  \ (3,2,2), \ (3,2,1), \  (3,1,1), \\
&&  (2,2,2), \  (2,2,1),  \  (2,1,1), \ (1,1,1).
\end{eqnarray*}
Of these, the only one that does not immediately force $I$ to have a non-trivial GCD in degree $t+3$ is $(3,2,1)$. In this case $a$ cannot be 2. If $a=3$, then $t+2 \geq a+1 = 4$, which provides $f_3 = f_4 = \cdots = f_{t+2} = 3$. If $a = 4$, then $t+2 \geq a+1 = 5$, which provides $f_4 = f_5 = \cdots = f_{t+2} = 3$. Hence in both cases we see that $I_{t+2}$ has a GCD of degree 3.

This shows that each of the above possibilities for $(f_{t+2}, f_{t+3}, f_{t+4})$ leads to the conclusion that  $I_{t+2}$ has a non-trivial GCD. Therefore Lemma \ref{lem-gcd} shows that $h$ must be a differentiable O-sequence up to degree $t+2$, a contradiction to the choice of $t+2$.
\medskip

 \noindent {\bf Case 2:} Assume $c_{e-t} \leq c_{e-t-1}$.  Then we have
 \[
 h_{e-t} - b_{e-t} \leq h_{e-t-1} - b_{e-t-1},
 \]
 so we obtain
 \[
 h_{e-t-1} - h_{e-t} \geq b_{e-t-1} - b_{e-t}.
 \]
 By symmetry, we get
 \begin{equation} \label{eqn1}
 h_{t+1} - h_t \geq b_{t+2} - b_{t+1}.
 \end{equation}
 Combining (\ref{h's}) with (\ref{eqn1}), we have
 \begin{equation}\label{eqn3}
 h_{t+2} - h_{t+1} > h_{t+1} - h_t \geq b_{t+2} - b_{t+1},
 \end{equation}
 which implies $ c_{t+2} > c_{t+1}$.  Hence by (\ref{c's}), we conclude
 \begin{equation}\label{eqn2}
 f_{t+2} > 0.
 \end{equation}
 Recall that the ``c'' line corresponds to the Hilbert function of $R/(I,L_1)$. We now claim that
\begin{equation} \label{eq-ccc}
c_{t+1} < c_{t+2} \leq c_{t+1} + 2.
\end{equation}
Indeed, this is clear by Macaulay's growth condition if $a=2$. Hence we may assume $3 \leq a \leq 4$, thus $t+2 \geq 4$. By Green's theorem,
 $h_4 \leq 33$ implies  that $c_4 \leq 13 = \binom{5}{4} + \binom{4}{3} + \binom{3}{2} + \binom{1}{1}$. Hence Macaulay's growth condition and  \cite[Lemma 4.13(b)]{N} provide $c_{i+1} - c_i \leq 3$ if $i \geq 4$ and equality is only possible if $(I,L_1)/(L_1)$ has a GCD of degree 3 in degree $i+1$. In case $a=3$, this conclusion is also true if $i = 3$. If $a=4$, then $t+2 \geq 5$. Hence we  get in either case that if $c_{t+2} = c_{t+1} + 3$, then $I_{t+2}$ has a GCD of degree 3. But then $h$ would not fail to be an SI-sequence in degree $t+2$, by Lemma \ref{lem-gcd}. This completes the proof of  Estimate (\ref{eq-ccc}).

We now consider the situation when $c_{t+2} = c_{t+1} + 2$. If this represents maximal growth of the Hilbert function, then we get that $(I,L_1)/(L_1)$ has a GCD of degree 2 in degree $t+2$, which, by Lemma \ref{lem-gcd}, implies that $h$ is a differentiable O-sequence up to degree $t+2$, contradicting again our assumption. Under what circumstances does $c_{t+2} = c_{t+1} + 2$ not represent maximal growth? Certainly we must have $a \geq 3$. Moreover, the identity
$3j = \binom{j+1}{j} + \binom{j}{j-1} + \binom{j-1}{j-2}$ implies that we also must have
\begin{equation} \label{eq-low}
c_j \geq 3j \quad {\rm whenever} \quad  j \leq t+1.
\end{equation}
Now, if $a=3$, this forces $(c_0,\ldots,c_{t+2}) = (1, 3, 6, 9, 12, \ldots, 3j, \ldots, 3 (t+1), 3 t + 5)$. This is a differentiable O-sequence. Since by our induction hypothesis the ``b'' line is an SI-sequence, we conclude that $h$ is a differentiable O-sequence up to degree $t+2$, contradicting our assumption.

If $a = 4$, then we get $(c_0,\ldots,c_{3}) = (1, 3, 6, 10)$ and, using also our assumption $h_4 \leq 33$,  always $c_j \leq 3j+1$. Together with Inequality (\ref{eq-low}), this implies
$$
3 \geq c_{t+1} - c_t \geq 2 = c_{t+2} - c_{t+1}.
$$
Since  the ``b'' line is an SI-sequence, a straightforward computation shows that the  growth from $h_{t+1} - h_t = (b_{t+1} - b_t) + (c_{t+1} - c_t)$  to $h_{t+2} - h_{t+1} = (b_{t+2} - b_{t+1}) + (c_{t+2} - c_{t+1})$ satisfies Macaulay's condition, which is again a contradiction to the choice of $t+2$.

Therefore,  we are left with considering the case where
 \[
 c_{t+2} = c_{t+1} +1.
 \]
This means that
 \[
 h_{t+2} - h_{t+1} = b_{t+2} - b_{t+1} + 1.
 \]
 Using (\ref{eqn3}) we get
 \[
 h_{t+1} - h_t = b_{t+2} - b_{t+1}.
 \]
 Combining these, we obtain
 \begin{equation} \label{eqn4}
 h_{t+2} - h_{t+1} = h_{t+1} - h_t + 1.
 \end{equation}
 But we have assumed that $h$ fails to be an SI-sequence in degree $t+2$, so (\ref{eqn4}) must represent a violation of maximal growth.  By Macaulay's theorem, this forces
 \begin{equation} \label{h diff}
 h_{t+1} - h_t \leq t+1 \ \ \hbox{ and } \ \ h_{t+2} - h_{t+1} \leq t+2.
 \end{equation}

Next we observe
\smallskip

\noindent
{\bf Claim:} $I_{t+1}$ does not have a non-trivial GCD.
\smallskip

Indeed, assume on the contrary that $I_{t+1}$ has a GCD, $F$, of degree $d > 0$. Using that $(I, F)_j = (F)_j$ if $j \leq t+1$ and also $t+1 \geq a \geq d > 0$, we get
\begin{eqnarray*}
\Delta h_{R/(I, F)} (t+1) = \Delta h_{R/(F)} (t+1) & = & \binom{t+3}{2} - \binom{t+3-d}{2} \\
& \geq & \binom{t+3}{2} - \binom{t+2}{2} \\
& = & t+2.
\end{eqnarray*}
Set $g_j = h_{R/I:F} (j)$. Recalling (\ref{h diff}), we obtain
 \[
 t+1 \geq h_{t+1} - h_t = g_{t+1-d} - g_{t-d} + \Delta h_{R/(F)} (t+1) \geq g_{t+1-d} - g_{t-d} + t+2,
 \]
 so
 \[
 g_{t+1-d} - g_{t-d} < 0.
 \]
 But $t+2 \leq \frac{e}{2}$ implies $t+1-d \leq \frac{e-d}{2}$. Since
 $R/(I:F)$ is Gorenstein with socle degree $e-d$, we must have $g_{t-d} \leq g_{t+1-d}$, a contradiction. This completes the proof of the claim.
\smallskip

We now consider the possibilities for the triple $(f_t, f_{t+1}, f_{t+2})$:
 \begin{itemize}
 \item[a.] $(*, 3, 3)$, $(*,2,2)$, $(*,1,1)$;

 \item[b.] $(3, 3, *)$, $(2, 2, *)$;

 \item[c.]  $(4, * , *)$;

 \item[d.] $(3,2,1)$.

 \end{itemize}

In subcase a., we have maximal growth on the ``$f$'' line from degree $t+1$ to degree $t+2$.   This means that  $((I,L_1, L_2)/(L_1, L_2))_{t+2}$ has a non-trivial GCD, hence so has $I_{t+1}$ (see Case 1 of the proof of Theorem \ref{unimodal thm}), which is impossible by the above Claim.

Similarly, in subcase b., maximal growth implies that $I_{t+1}$ has a non-trivial GCD,  a contradiction.

The condition in subcase c.\ is equivalent to $a=4$ and $t=3$. Thus, (\ref{h diff}) provides
\[
h_4 \leq h_3 + 4 = 24.
\]
Using Green's theorem twice, we get that $f_4 \leq 1$. But $f_5 = f_{t+2} > 0$, so $f_4 = f_5 = 1$. It follows that $I_5$ has a GCD of degree 1, and we conclude as in subcase a.

 Finally, we turn to subcase d. We have $a \in \{3, 4\}$.
Assume first that $a=3$.
Then by Proposition \ref{prop-key}, either $I$ has a GCD of degree 2 in degree $t+1$ or else $I $ has up to degree $t+1$ exactly two minimal generators. This means that either $I_{t+1}$ has a non-trivial GCD or $I$ is up to degree $t+1$ generated by a regular sequence of length two.  The former case contradicts the above Claim. So the Hilbert function of $R/I$ agrees with that of a complete intersection of height 2 through degree $t+1$, and $I$ possibly picks up new generators in degree $t+2$.  But then $h$ must be below the Hilbert function of the complete intersection, contradicting the assumption that $h_{t+2}$ was the value that violated the SI condition.

Second, assume $a = 4$. If $t \geq 5$, then $I_t$ has a GCD of degree 3. If $t = 4$, then applying Proposition \ref{prop-key} we see that either $I_t$ has a GCD of degree 3 or $I$ has up to degree $t = 4$ exactly two minimal generators of degree 4.

Let us start with considering the case where $I_t$ has a  GCD, $F$, of degree 3. Then using the arguments in case $a=4$ of the proof of Theorem \ref{unimodal thm}, we conclude that $(I,F)_{t+1}$ has a GCD, $D$, of degree 2 or $(I, F)$ has only two minimal generators up to degree $t+1$. In the former case, $D$ is a common divisor of $I_{t+1}$, and we get a contradiction to the above claim. If in the latter case $(I, F)$ has a non-trivial GCD, we argue as above. Otherwise $(I, F)$ is up to degree $t+1$ a complete intersection $J$ of height two that is generated by forms of degree 2 and $t+1$. Hence, an easy computation yields
\[
\Delta h_{R/(I, F)} (t+1) =  \Delta h_{R/J} (t+1) = 2 t + 2 > t+1.
\]
This provides a contradiction to (\ref{h diff}) as in the proof of the above claim.

It remains to discuss the case where $a = t= 4$ and $I$ has exactly two minimal generators of degree 4. This implies $c_4 = 13$ and $h_4 = 33$. Using (\ref{h diff}) we get $h_5 \leq 38$, thus $c_5 \leq 11$ by Green's theorem. But $f_5 =2$ and Green's theorem provide $c_5 \geq 11$. Hence we obtain
\[
c_5 = 11 < c_4 = 13.
\]
Since $t+2 = 6 \leq \frac{e}{2}$, we get $e \geq 12$. Thus, the arguments at the beginning of Case 2 imply that $f_{j+2} > 0$ for some $j \geq 7$. Since $f_6 = 1$, we get in particular $f_6 = f_7 = 1$, which implies that $I_{7}$ has a GCD of degree 1 (see Case 1 of the proof of Theorem \ref{unimodal thm}), and so does  $I_5$; but this contradicts the above Claim.

The proof is now complete.
\end{proof}

At this point the reader might wonder whether it is possible to extend Theorem \ref{more general} to guarantee the SI-property. However, the methods of the previous proof are not enough to establish such an extension. Using the notation as in Theorem \ref{more general} and the above proof, its arguments in Case 1 and Case 2 extend to the more general situation, except when we have in Case 2 that $t+2 < s$. In this situation a new argument  seems to be  needed.

\bigskip

\noindent {\bf Acknowledgements:}  The authors express their sincere thanks to Tony Iarrobino for his extremely thorough reading of an earlier version of this work, which led to the discovery of a serious error and to a substantial revision, resulting in the present paper.  The research performed in this paper was started during the visit of the
first and third authors to the second one at the University of Kentucky in
spring 2006.

The authors also thank the referee for helpful comments.

\end{document}